\documentclass{amsart}
\usepackage{amsmath,amssymb,graphicx,bbm}

\newcommand{\mathe}{\mathrm{e}}

\newcommand{\tmop}[1]{\ensuremath{\operatorname{#1}}}
\newcommand{\tmstrong}[1]{\textbf{#1}}
\newcommand{\tmtextit}[1]{{\itshape{#1}}}

\begin{document}

\title[Mean Staircase of the Riemann Zeros]{Mean Staircase of the Riemann Zeros: a comment on the Lambert W
function and an algebraic aspect}
\author{Davide a Marca}
\address{D. a Marca, CERFIM, Research Center for Mathematics
and Physics, PO Box 1132, 6600 Locarno, Switzerland}
\email{damarca@ticino.com}
\author{Stefano Beltraminelli}
\address{S. Beltraminelli, CERFIM, Research Center for Mathematics
and Physics, PO Box 1132, 6600 Locarno, Switzerland}
\email{stefano.beltraminelli@ti.ch}
\author{Danilo Merlini}
\address{D. Merlini, CERFIM, Research Center for Mathematics and
Physics, PO Box 1132, 6600 Locarno, Switzerland}
\email{merlini@cerfim.ch}
\label{I1}\label{I2}\label{I3}
\date{17.1.2009}
\subjclass{11M26}
\keywords{Riemann zeta function, LambertW function, Riemann zeros, harmonic oscillator, Riemann Hypothesis}
\begin{abstract}
  In this note we discuss explicitly the structure of two simple set of zeros
  which are associated with the mean staircase emerging from the zeta function
  and we specify a solution using the Lambert W function. The argument of it
  may then be set equal to a special $N \times N$ classical matrix (for every
  $N$) related to the Hamiltonian of the Mehta-Dyson model. In this way we
  specify a function of an hermitean operator whose eigenvalues are the
  ``trivial zeros'' on the critical line. The first set of trivial zeros is
  defined by the relations $\tmop{Im} \left( \zeta \left( \frac{1}{2} + i
  \cdot t \right) \right) = 0 \wedge \tmop{Re} \left( \zeta \left( \frac{1}{2}
  + i \cdot t \right) \right) \neq 0$ and viceversa for the second set. (To
  distinguish from the usual trivial zeros $s = \rho + i \cdot t = - 2 n$, $n
  \geqslant 1$ integer)
\end{abstract}
\maketitle
\dedicatory{\itshape{This (heuristic, non rigorous) research note is dedicated to the
  international Swiss-Italian mathematician and physicist Professor Dr. Sergio
  Albeverio on the occasion of his seventieth birthday; a friend and for years
  the scientific director of Cerfim (Research Center for Mathematics and
  Physics of Locarno), situated opposite the "Rivellino"}{\footnote{The Bastion
  "Il Rivellino", situated opposite Cerfim, is 99.9\% attribuable to Leonardo
  da Vinci (1507).}}.}
  
\vskip 1cm

\section{Introduction: a search for an hermitean operator associated with the
Riemann Zeta Function}

There is much interest in understanding the complexity related to the Riemann
Hypothesis and concerned with the location and the structure of the non
trivial zeros of the Riemann zeta function $\zeta (s)$ where $s = \rho + i
\cdot t$ is the complex variable. Following a suggestion of Hilbert and Polya,
in recent years many efforts have been devoted to a possible construction of
an hermitean operator having as eigenvalues the imaginary parts $t_n$ of the
non trivial zeros of $\zeta$ ($\zeta$ being meromorphic, the zeros are
countable). These are given by the solutions of the equation $\zeta \left(
\rho_n + i \cdot t_n \right) = 0$, $n = 1, 2, \ldots,$. If $\rho_n =
\frac{1}{2}$ for all $n$, then all the zeros lie on the critical line (the
Riemann Hypothesis is true); the program is then to find an hermitean
``operator'' $T$ such that $T \cdot \varphi_n = t_n \cdot \varphi_n$ in some
appropriate space ($\varphi_n$ would be the $n^{\tmop{th}}$ eigenvector of
$T$). There are today many strategies in the direction of constructing such an
operator and in the sequel we will comment on some (among many others) very
stimulating works on the subject. In {\cite{pitkanen-2003}}, Pitk\"anen's
heuristic work goes in the direction of constructing orthogonality relations
between eigenfunctions of a non hermitean operator related to the
superconformal symmetries; a different operator than the one just mentioned
has also been proposed in {\cite{castro-2001}} by Castro, Granik and Mahecha
in terms of the Jacobi Teta series and an orthogonal relation among its
eigenfunctions has also been found. In the rigorous work by Elizalde et.al
{\cite{elizalde-2003}} some problems with those approches have been pointed
out. In a work of some years ago Julia {\cite{julia-1989}} proposed a
fermionic version of the zeta function which should be related to the
partition function of a system of $p$-adic oscillators in thermal equilibrium.
In two others pioneering works of these years, Berry and Keating
{\cite{berry-1999a,berry-1999b}} proposed an interesting heuristic operator to
study the energy levels $t_n$ (the imaginary parts of the non trivial zeros of
the zeta function). The proposed Hamiltonian has a very simple form given, on
a dense domain, by: $H = p \cdot x + \frac{1}{2}$, where
\begin{equation}
  \label{pBerry} p = \left( \frac{1}{i} \right)  \frac{\partial}{\partial x}
\end{equation}
in one dimension. As explained by the authors, the difficulty is then to
define appropriate spaces and boundary conditions to properly determine $p$
and $H$ as hermitean operators. In such an approach the heuristic appearance
of ``instantons'' is also discussed. In another important work Bump et al.
{\cite{bump-2000}} introduced a local Riemann Hypothesis and proved in
particular that the Mellin transform of the Hermite polynomials (associated
with the usual quantum mechanical harmonic oscillator) contain as a factor a
polynomial $p_n (s)$, corresponding to the $n$-energy eigenstate of the
oscillator, whose zeros are exactly located on the critical line $\sigma =
\frac{1}{2}$. The relation of the polynomials $p_n (s)$ with some truncated
approximation of the entire funcion $\xi (s)$ (the Xi function), related to
the Riemann zeta function seems to be still lacking. Others important
mathematical results concerning the non trivial Riemann zeros, have been
obtained by many leading specialist (see among others the work by Connes
{\cite{connes-1999}} and the work by Albeverio and Cebulla
{\cite{albeverio-2007}}). For a recent work on the xp hamiltonian see G.
Sierra {\cite{sierra-2007}}.

Let us also mention that for the nontrivial zeros of zeta an interesting
equation has been proposed originally by Berry and Keating in
{\cite{berry-1999a}}. In fact, remembering the definition of $\xi (s)$ where
$\zeta (s)$ is the Riemann zeta function, given by:
\begin{equation}
  \label{Xi} \xi (s) = \frac{1}{2} \cdot s \cdot \left( s - 1 \right) \cdot
  \pi^{^{- \frac{s}{2}}} \cdot \Gamma \left( \frac{s}{2} \right) \cdot \zeta
  \left( s \right)
\end{equation}
the equation for possible zeros of $\xi$ proposed in {\cite{berry-1999a}} is
given by:
\begin{equation}
  \label{gs} \frac{\pi^{\frac{s}{2}}}{\Gamma \left( \frac{s}{2} \right)} +
  \frac{\pi^{\frac{1 - s}{2}}}{\Gamma \left( \frac{1 - s}{2} \right)} = 0
\end{equation}
As stated by the authors, Eq. (\ref{gs}) could be considered as a
``quantization condition''. Unfortunately, as mentioned in
{\cite{berry-1999a}}, Eq. (\ref{gs}) possesses complex zeros and so can not be
used to provide an hermitean operator which would generate the non trivial
zeros of $\zeta$. The content of our note is concerned with the ``mean
staircase'' of the Riemann zeros: we first specify the two set of trivial
zeros on the critical line related to it and point out an explicit
construction using the Lambert W function; then we introduce a specific
argument (a $n \times n$ hermitean matrix $H$, describing a discrete harmonic
oscillator with creation and annihilation ``operators'' $a$ and $a^{\ast}$
such that $[a, a^{\ast}] = - 2$) into the Lambert W function. We obtain then,
for the trivial zeros, the goal that the ``Polya-Hilbert program'' has for the
non trivial zeros.

\section{The mean staircaise of the Riemann zeros and the trivial zeros on
the critical line associated with it}

Let $\xi (s)$ be the xi function given by (\ref{Xi}). If $N (t)$ denotes the
number of zeros of $\xi$ in the critical strip of height smaller or equal to
$t$, and if $S (t) \equiv \frac{1}{\pi} \arg \left( \zeta \left( \frac{1}{2} +
i t \right) \right)$, then {\cite{titchmarsh-1986}}
\begin{equation}
  \left. \label{Nt} N (t) = \langle N (t) \right\rangle + S \left( t \right) +
  O \left( \frac{1}{t} \right),
\end{equation}
where
\begin{equation}
  \label{mNt} \left\langle N (t) \right\rangle = \frac{t}{2 \pi} \cdot \left(
  \ln \left( \frac{t}{2 \pi} \right) - 1 \right) + \frac{7}{8}
\end{equation}
$\left\langle N (t) \right\rangle$, the ``bulk contribution'' to $N$, is
called the ``mean staircase of the zeros'' (cfr. {\cite{titchmarsh-1986}}).
The fluctuations of the number of zeros around the mean staircase, are given
by the function $S (t)$. It is known {\cite{titchmarsh-1986}} that $S (t) = O
(\ln t)$ without assuming RH while, assuming RH is true, it is known that $S
(t) = O \left( \frac{\ln t}{\ln \left( \ln t \right)} \right)$. At this point,
since our remark has mainly to do with $\left\langle N (t) \right\rangle$, we
will set $S (t) = 0$ in Eq.(\ref{Nt}). We shall study the relation $N (t) =
\left\langle N (t) \right\rangle + O \left( \frac{1}{t} \right)$. The two sets
(which we call here ``trivial zeros on the critical line'') of interest are
defined by the above mean staircase as follows. The first set is given by the
zeros of $\tmop{Im} \zeta \left( \frac{1}{2} + i t \right)$ alone i.e. such
that $\tmop{Re} \zeta \left( \frac{1}{2} + i t \right) \neq 0$ (the first set
of trivial zeros on the critical line). The second set is given by the zeros
of $\tmop{Re} \zeta \left( \frac{1}{2} + i t \right)$, such that $\tmop{Im}
\zeta \left( \frac{1}{2} + i t \right) \neq 0$. For the first set:
\begin{equation}
  \tmop{Im} \zeta \left( \frac{1}{2} + i t_n^{\ast} \right) = 0 \wedge
  \tmop{Re} \zeta \left( \frac{1}{2} + i t_n^{\ast} \right) \neq 0
\end{equation}
Then since $\tmop{Im} \zeta \left( \frac{1}{2} + i t_n^{\ast} \right) = 0$ we
have that $\pi^{- 1} \arg \zeta \left( \frac{1}{2} + i t_n^{\ast} \right) =
\pi^{- 1}  \left( - \pi n \right)$ and $\pi^{- 1} \arg \xi = 0$ is given by
those $t_n^{\ast}$ such that
\begin{equation}
  N (t_n^{\ast}) \cong \left\langle N (t_n^{\ast}) \right\rangle =
  \frac{t_n^{\ast}}{2 \pi} \cdot \left( \ln \left( \frac{t_n^{\ast}}{2 \pi}
  \right) - 1 \right) + \frac{7}{8} = n, n \tmop{integer} \geqslant 1
\end{equation}
at large values of $t$ or $n$. So the nonlinear equations to be solved which
should give the values where only $\tmop{Im} \zeta \left( \frac{1}{2} + i t
\right)$ vanishes, i.e $\{t_n^{\ast} \}$, is given by:
\begin{equation}
  \label{mNts} \left\langle N (t_n^{\ast}) \right\rangle = n, n \tmop{integer}
  \geqslant 1
\end{equation}
while for the second set
\begin{equation}
  \label{mNtss} \left\langle N (t_n^{\ast \ast}) \right\rangle = n -
  \frac{1}{2}, n \tmop{integer} \geqslant 1
\end{equation}
(The first set has been known for a long time and constitutes the Gram points,
$\sin (\theta) = 0$, where $\theta$ is the phase of the $\zeta$ function,
while for the second set one has $\cos (\theta) = 0$).

We note, the values of interest are given by the abscissa of the intersection
points between the staircase (Eq.(\ref{mNt})) and the two functions $\pi^{- 1}
\arg \xi \left( \frac{1}{2} + i t \right)$ and $\text{$\pi^{- 1} \arg \xi
\left( \frac{1}{2} + i t \right)$} - \frac{1}{2}$ . The plot of Fig 1
illustrates the situation for some low lying zeros. The values for
$t_n^{\ast}$ lie mostly in between the exact value of the Riemann zeros $t_{n
- 1}$ and $t_n$, but it is known that the Gram law fails for the first time at
$t = 282.4$ (``first istanton'' according to {\cite{berry-1999a}}). The
solution of the above equation which gives $t_n^{\ast}$, $t_n^{\ast \ast}$
using a very special function (the LambertW function, see
{\cite{corless-1996}}) is given below.

\begin{figure}[h]
  \includegraphics[width=10cm]{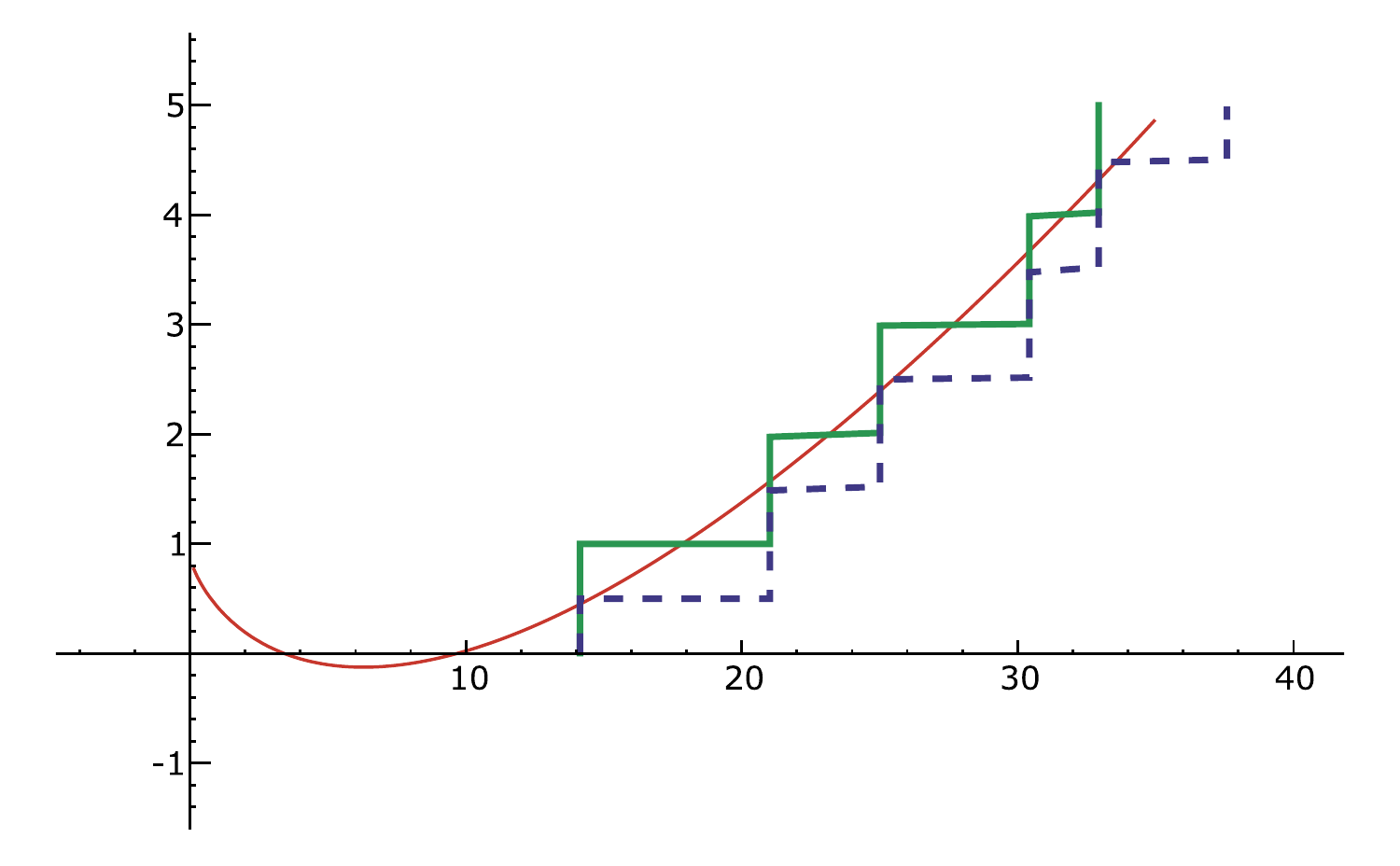}
  \caption{The plot of $\left\langle N (t) \right\rangle$ (continuos curve),
  of $N (t)$ (full stair) and $\text{$N (t)$} - \frac{1}{2}$ (intermittent
  stair)}
\end{figure}

\section{An exact solution for the sequence $t_n^{\ast}$ and $t_n^{\ast
\ast}$}

The equation corresponding to (\ref{mNts}), may be written in the form
\begin{equation}
  \label{lamba} \left( \frac{t}{2 \pi \mathe} \right)^{\frac{t}{2 \pi \mathe}}
  = \mathe^{\frac{n - \frac{7}{8}}{\mathe}}
\end{equation}
and the equation corresponding to (\ref{mNtss}) in the form
\begin{equation}
  \label{lambb} \left( \frac{t}{2 \pi \mathe} \right)^{\frac{t}{2 \pi \mathe}}
  = \mathe^{\frac{n - \frac{1}{2} - \frac{7}{8}}{\mathe}}
\end{equation}
so that introducing the new variables $x = \exp \left( \frac{n -
\frac{7}{8}}{\mathe} \right)$ resp. $x = \exp \left( \frac{n - \frac{1}{2} -
\frac{7}{8}}{\mathe} \right)$ we obtain the equation (from (\ref{lamba}) and
(\ref{lambb}), $x > 0$)
\begin{equation}
  W (x) \cdot \exp (W (x)) = x
\end{equation}
The function $W (x)$ is called the Lambert W function and has been studied
extensively in these recent years. In fact such an equation appears in many
fields of science. In particular the use of such a function has appeared in
the study of the wave equation in the double-well Dirac delta function model
or in the solution of a jet fuel problem. See {\cite{corless-1996}} for an
important work on the subject. Moreover the Lambert W function appears also in
combinatorics as the generating function of trees and as explained in
{\cite{corless-1996}} the W function has many applications, even if the
presence of the W function often goes unrecognized.

The Lambert W function has many complex branches; of interest here is the
principal branch of W which is analytic at $x = 0$. So, the solution of
(\ref{mNts}, \ref{mNtss}) is given by
\begin{equation}
  \label{ts} t_n^{\ast} = 2 \pi \mathe \cdot \exp \left( W \left( \frac{n -
  \frac{7}{8}}{\mathe} \right) \right)
\end{equation}
\begin{equation}
  \label{tss} t_n^{\ast \ast} = 2 \pi \mathe \cdot \exp \left( W \left(
  \frac{n - \frac{1}{2} - \frac{7}{8}}{\mathe} \right) \right)
\end{equation}
We have thus specified, with the help of the LambertW function, the sequences
$\{t_n^{\ast} \}$ resp. $\{t_n^{\ast \ast} \}$, which are the zeros of \
$\tmop{Im} \left( \zeta \left( \frac{1}{2} + i \cdot t \right) \right)$ such
that $\tmop{Re} \left( \zeta \left( \frac{1}{2} + i \cdot t \right) \right)
\neq 0$ resp. $\tmop{Re} \left( \zeta \left( \frac{1}{2} + i \cdot t \right)
\right) = 0$ and $\tmop{Im} \left( \zeta \left( \frac{1}{2} + i \cdot t
\right) \right) \neq 0$.

It should be noted here that in Eq.(\ref{mNts}), $n$, which would correspond
to the exact value of a true zero value $t_n$ (non trivial zero) of the
$\zeta$ function would not be an integer $n$ or $n - \frac{1}{2}$ since we
have replaced in Eq. (\ref{mNt}) $S (t)$ by zero. For the first few low zeros
(the true zeros), it may be observed numerically that the corresponding
values, let say $n^{\ast}$, are randomly distribued mostly between two
consecutive integers, but the mean values are nearby the integers plus
$\frac{1}{2}$. A calculation with some zeros gives a mean value of 0.49
instead of 0.5. So, in average it seems that the behavior of the true zeroes
$t_n$ ``follows'' more the pattern of the set $t_n^{\ast \ast}$. In the
similar way the zeros of the first set, i.e. $t_n^{\ast}$, lie mostly in
between two non trivial zeros of $\zeta$ but of course it is known that there
are very complicated phenomena associated with the chaotic behavior of the non
trivial zeros of the Riemann $\zeta$ function.

As an example, the first of the istantons, corresponding to $n = 126$, cited
above, is located at the value of $t = 282.4 \ldots$. On the table below we
give the values of $t_n^{\ast}$ and of $t_n$ of a true zero around $t = 280$.

\begin{table}[h]

  \begin{center}
    \begin{tabular}{l}
      $t_{126} = 279.22925$\\
      $t_{126}^{\ast} = 280.80246$\\
      $t_{127}^{\ast} = 282.4547596$\\
      $t_{127} = 282.4651147$\\
      $t_{128} = 283.211185$\\
      $t_{128}^{\ast} = 284.1045158$\\
      $t_{129} = 284.8359639$
    \end{tabular}
  \end{center}
  \caption{}
\end{table}

From those numerical computations we see that two consecutive zeros of
$\tmop{Im} (\zeta)$ alone are followed by two consecutive true zeros, that is
$t_{127}^{\ast}$ anticipates $t_{127}$. The difference between the two
subsequent $t$ values is very small and given by $\Delta t = 0.0103$. The
phase change is given by $i \pi$ as illustrated on the plot of $\tmop{Im}
\left( \ln \left( \zeta \left( \frac{1}{2} + i t \right) \right) \right)$
(step curve) and that of $\tmop{Im} \left( \zeta \left( \frac{1}{2} + i t
\right) \right)$.

\begin{figure}[h]
  \includegraphics[width=10cm]{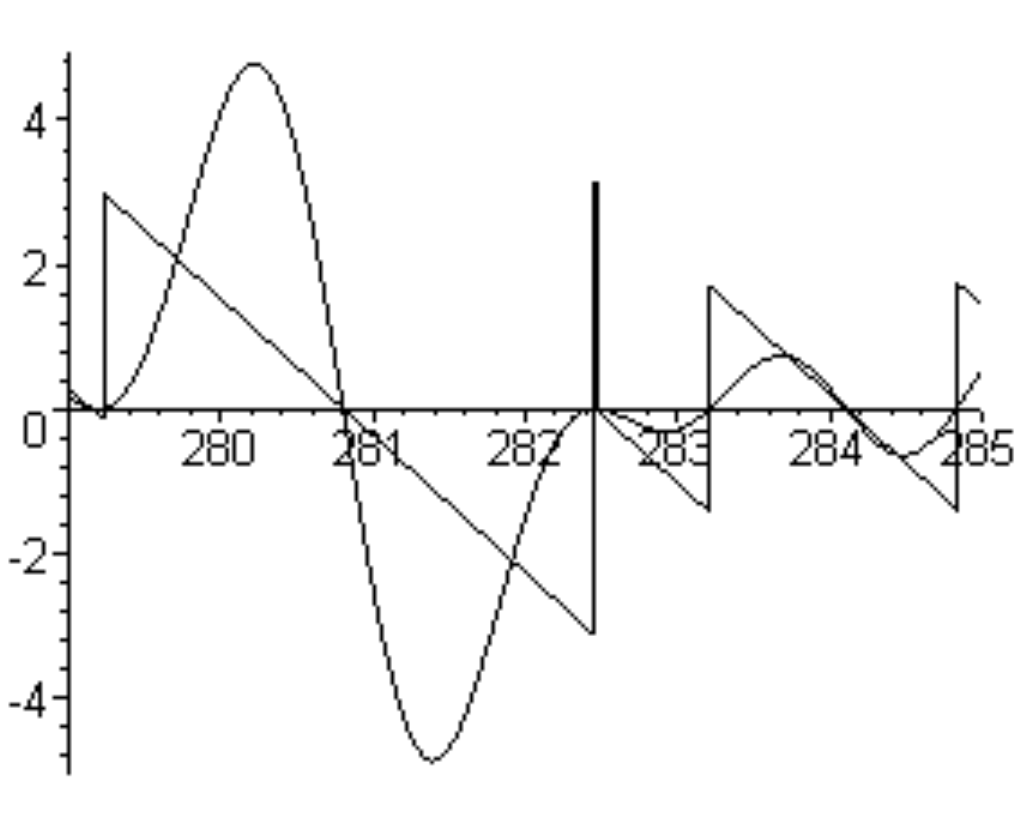}
  \caption{The first istanton}
\end{figure}

For the first 500 energy levels, that is for values of $t$ from 0 to $t =
811.184 \ldots$ (level number $n = 500$), it may be seen that there are 13
istantons (in the language of {\cite{berry-1999a}}), all with a Maslov phase
change of $+ i \pi$ or of $- i \pi$. The width is usally small but it is
larger for the istanton located at $t = 650.66$ ($n$ corresponding to $379$),
where this time $\Delta t = 0.31 \ldots$. Returning now to the trivial zeros
(the two sets $t^{\ast}_n$, $t^{\ast \ast}_n$ defined above), we note the
elementar relation which follow from (\ref{mNts}) and (\ref{mNtss}), and given
by:
\begin{equation}
  \label{av1} \frac{t_n^{\ast} + t_{n + 1}^{\ast}}{2} = t_{n +
  \frac{1}{2}}^{\ast \ast}
\end{equation}
and
\begin{equation}
  \label{av2} \frac{t_{n - \frac{1}{2}}^{\ast \ast} + t_{n +
  \frac{1}{2}}^{\ast \ast}}{2} = t_n^{\ast}
\end{equation}
Eq.(\ref{av1}, \ref{av2}) say that the zeros of the real part alone are
obtained by those of the imaginary part alone by simple average and viceversa.
The two sequences are regularly spaced and the mean distance between two
trivial zeros at the height $t$, as the mean staircase indicates
(Eq.(\ref{mNt})), is given approximatively by:
\begin{equation}
  \frac{t}{\left\langle N (t) \right\rangle} = \frac{2 \pi}{\log \left(
  \frac{t}{2 \pi} \right)} = \frac{2 \pi}{\log (n)}
\end{equation}
for $t$ and $n$ large

Before proposing an hermitean operator for the sequences of the trivial zeros
it is important to investigate a possible ``quantization condition'' for the
non trivial zeros. For this we start with the Riemann symmetry of the $\zeta$
function.

From the exact relation for the $\xi$-function given by:
\begin{align}
  \xi (s) &= \frac{1}{2} \pi^{- \frac{s}{2}} \Gamma \left( \frac{s}{2} \right)
  \zeta (s) s (s - 1) \\
&= \xi (1 - s) = \frac{1}{2} \pi^{- \frac{1 - s}{2}}
  \Gamma \left( \frac{1 - s}{2} \right) \zeta (1 - s) (1 - s) (1 - s - 1) \notag
\end{align}
$s \in \mathbbm{C}$, we have that
\begin{equation}
  \label{ZG} \pi^{- \frac{s}{2}} \Gamma \left( \frac{s}{2} \right) \zeta (s) =
  \pi^{- \frac{1 - s}{2}} \Gamma \left( \frac{1 - s}{2} \right) \zeta (1 - s)
\end{equation}
In equation (\ref{ZG}) we limit ourselves to consider the values $s = \rho + i
t = \frac{1}{2} \downarrow + i t$, $t \in \mathbbm{R}$, and thus $1 - s =
\frac{1}{2} \uparrow - i t$; moreover we are interested in high values of $t$
so that we may use the Stirling's formula for the Gamma function given by:
\begin{equation}
  \label{Ge} \Gamma (x) \cong \left( 2 \pi \right)^{\frac{1}{2}} x^{x -
  \frac{1}{2}} \mathe^{- x}
\end{equation}
as $x \rightarrow \infty$. From (\ref{ZG}) and (\ref{Ge}) we then obtain
(asymptotically for $t \rightarrow \infty$)
\begin{align}
 &\exp \left( i \pi \left( \left( \frac{t}{2 \pi} \right)  \left( \ln \left(
  \frac{t}{2 \pi} \right) - 1 \right) - \frac{1}{8} \right) + i \arg \left(
  \zeta \left( \frac{1}{2^+} + i t \right) \right) \right) =\\
&\exp \left( - i \pi \left( \left( \frac{t}{2 \pi} \right)  \left( \ln \left(
  \frac{t}{2 \pi} \right) - 1 \right) - \frac{1}{8} \right) + i \arg \left(
  \zeta \left( \frac{1}{2^-} - i t \right) \right) \right) \notag
\end{align}
Since 
\begin{align}
\exp \left( i \arg \left( \zeta \left( \frac{1}{2^{}} \uparrow - i t
\right) \right) \right) &= \exp \left( i \arg \left( \zeta \left(
\frac{1}{2^{}} \downarrow + i t \right) + i \pi \right) \right) \notag \\
&= - \exp
\left( i \arg \left( \zeta \left( \frac{1}{2^{}} \downarrow + i t \right)
\right) \right) \notag
\end{align}
we then have, taking the limit $\rho = \frac{1}{2}
\downarrow = \frac{1}{2}$, that:
\begin{equation}
  \label{cospsi} \cos (\Psi) = 0 \tmop{where} \Psi = \frac{t}{2}  \left( \ln
  \left( \frac{t}{2 \pi} \right) - 1 \right) - \frac{\pi}{8} + \arg \left(
  \zeta \left( \frac{1}{2} + i t \right) \right)
\end{equation}
Thus $\Psi = \pi \left( n + \frac{1}{2} \right)$. We then obtain:
\[ \frac{t}{2 \pi}  \left( \ln \left( \frac{t}{2 \pi} \right) - 1 \right) -
   \frac{1}{8} + \frac{1}{\pi} \arg \left( \zeta \left( \frac{1}{2} + i t
   \right) \right) = n - \frac{1}{2} \]
hence
\begin{equation}
  \frac{t}{2 \pi}  \left( \ln \left( \frac{t}{2 \pi} \right) - 1 \right) +
  \frac{7}{8} + \frac{1}{\pi} \arg \left( \zeta \left( \frac{1}{2} + i t
  \right) \right) = n + \frac{1}{2}
\end{equation}
Eq.(\ref{cospsi}) may be seen as an approximate ``quantum condition'' for the
true Riemann zeros, but it is only a consequence of the Riemann symmetry (Eq.
\ref{ZG}). In fact, if in equation (\ref{cospsi}) we neglect the last term
$\arg (\zeta)$, then (\ref{cospsi}) has as a solution the second set of
trivial zeros $\{t_n^{\ast \ast} \}$. It is true, as remarked by Berry and
Keating, that their Eq.(\ref{gs}) has complex zeros wich are not the Riemann
zeros, but it should be remarked that if in Eq.(\ref{gs}) we set $\tmop{Re}
(s) = \frac{1}{2}$ then Eq.(\ref{gs}) reduces to Eq.(\ref{cospsi}) without the
fluctuation term $\arg (\zeta)$; so the solution of Berry and Keating
Eq.(\ref{gs}) for $\tmop{Re} (s) = \frac{1}{2}$ is the same as the second set
of trivial zeros $\{t_n^{\ast \ast} \}$ we have specified.

Below the plots of the left hand side of Eq.(\ref{cospsi}), with and without
the term $\arg \zeta \left( \frac{1}{2} + i t \right)$. As an illustration, we
may observe on the plot the first istanton discussed above and the second one.
In fact the maximum of the function which gives $t^{\ast}$ (Eq. \ref{cospsi}
without the term $\arg \zeta$) is outside the plot of the step function given
by (\ref{cospsi}) (the true function). This is visible on the plot near $t =
282$ and near $t = 296$ (the second istanton). This conclude our remark on
(\ref{gs}) and Eq.(\ref{cospsi}). In the next section, we shall costruct an
hermitean operator whose eigenvalues are the trivial zeros of the zeta
function on the critical line.

\begin{figure}[h]
  \includegraphics[width=10cm]{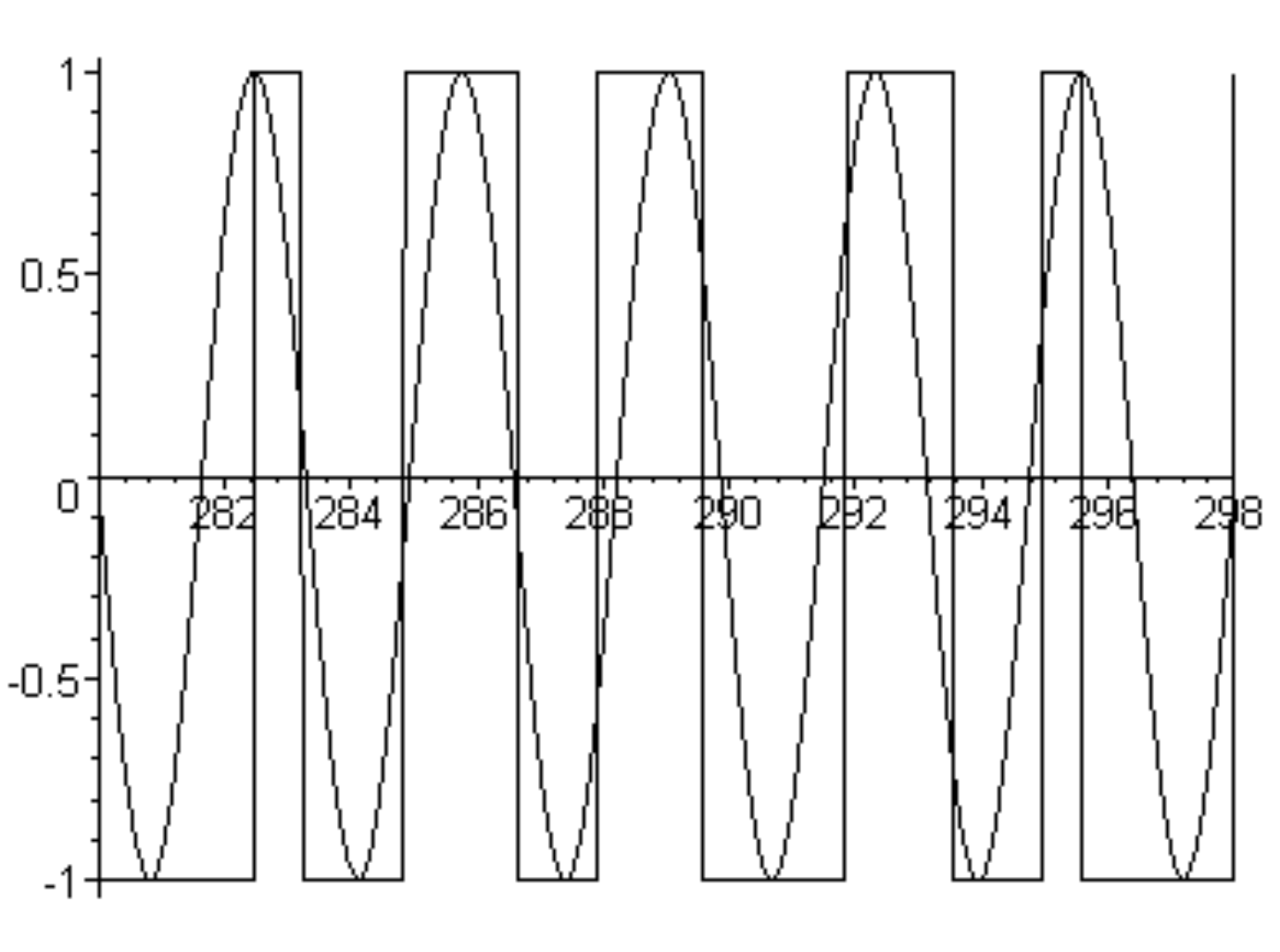}
  \caption{Plot of the function $\cos (\Psi)$ of Eq.(\ref{cospsi}) with the
  term $\arg \left( \zeta \left( \frac{1}{2} + i t \right) \right)$ (step
  function) and without that term.}
\end{figure}

Now in (\ref{ts}) and (\ref{tss}) the value of a trivial zero ($t_n^{\ast}$ or
$t_n^{\ast \ast}$) is given through his index $n$ by means of the Lambert W
function so that such zeros are related in a non linear way to the integers
$n$, i.e. in principle to the spectrum of an harmonic oscillator. So, for the
trivial zeros, no boundary condition is needed here, since they are obtained
by means of (\ref{ts}) and (\ref{tss}) in the large $t$ limit. At this moment
we are free to introduce a hermitean matrix which may generates the trivial
zeros.

\section{An Hermitean operator (matrix) associated with the mean staircase
(trivial zeros) of the Riemann Zeta function}

As remarked above, in (\ref{ts}) and (\ref{tss}) the only ``quantal number''
is the index $n$ of the trivial zeros and the construction may be given using
a hermitean $n \times n$ matrix $H$, for any $n$, at our disposal and related
to the classical one dimentional many body system whose fluctuation spectrum
around the equilibrium positions is that of the harmonic oscillator. In fact,
the one dimensional Mehta-Dyson model of random matrices (which may be seen as
a classical Coulomb system with $n$ particles) has, at low temperature an
energy fluctuation spectrum given by the integers and it is possible to
introduce classical annihilation and creation operators, as studied in
{\cite{merlini-1999}} (a short discussion is presented in the Appendix). The
matrix elements of the associated hermitean matrix are then functions of the
zeros of the Hermite polynomials; in this case we do not have a Hilbert space
and no Schr\"odinger Equation will be associated with the Lambert W function.
Another direction, i.e. that of introducing a Schr\"odinger Equation to
describe the trivial zeros may in principle be obtained as an application of
the results given by G. Nash {\cite{nash-1985}}; this because for large $n$,
as it is known, (\ref{ts}) and (\ref{tss}) give the behavior
({\cite{titchmarsh-1986}}, pag. 214) related to the asymptotic behavior of the
LambertW function:
\begin{equation}
  t_n = \frac{2 \pi n}{\ln (n)}, n \rightarrow \infty
\end{equation}
and thus the spectrum appears in fact as a one where the associated
Schr\"odinger Equation contains a Gaussian type of potential
{\cite{nash-1985}}. Here we will consider the matrix formulation: the point
may seem to be somewhat artificial but the hermitean matrix we will use
(specified in the Appendix) is related to the Mehta-Dyson model, the
``starting point'' of the random matrix theory. To do this, we begin to write
(\ref{lamba}) in a slighly different form using the Stirling formula for the
Gamma function of real argument given by:
\[ \Gamma (x) = \left( 2 \pi \right)^{\frac{1}{2}} x^{x - \frac{1}{2}}
   \mathe^{- x} \tmop{as} x \rightarrow \infty \]

We then have that, as $t \rightarrow \infty$,
\begin{align}
  \ln \left( \Gamma \left( \frac{t}{2 \pi \mathe} + \frac{1}{2} \right)
  \right) &= \frac{t}{2 \pi \mathe} \ln \left( \frac{t}{2 \pi \mathe} - 1
  \right) + \frac{7}{8} + \frac{1}{2} \ln \left( 2 \pi \right) - \frac{7}{8} \\ &=
  n^{\ast} + \frac{1}{2} \ln \left( 2 \pi \right) - \frac{7}{8} = n^{\ast} +
  \theta \notag
\end{align}
where $\theta = \frac{1}{2} \ln \left( 2 \pi \right) - \frac{7}{8}$.

Thus introducing the operator $T = T (H)$ whose eigenvalues should be the
trivial energy levels (for the first as well as for the second set defined by
(\ref{ts}) and (\ref{tss}) as well as $H$, the hermitean matrix given in the
Appendix and related to the Mehta-Dyson model, we may write the following
heuristic matrix equation:
\begin{equation}
  \label{GT} \Gamma \left( \frac{T}{2 \pi} + \frac{I}{2} \right) = \mathe^{H +
  \theta}
\end{equation}
where $I$ is the unit matrix. Eq.(\ref{GT}) is the equation for $T$, giving
the trivial zeros. The inversion of this formula (if it is possible to take
it) yields heuristically:
\begin{equation}
  T = T (H) = 2 \pi \left( \Gamma^{- 1} \left( \mathe^{H + \theta} \right) -
  \frac{I}{2} \right)
\end{equation}
To conclude, if $H \varphi_n = \left( n + \frac{1}{2} \right) \varphi_n$,
where $\varphi_n$ is the $n^{\tmop{th}}$ eigenfunction of $H$, then
\begin{equation}
  T \varphi_n = 2 \pi \left( \Gamma^{- 1} \left( \mathe^{H + \theta} \right) -
  \frac{I}{2} \right) \varphi_n = 2 \pi \left( \Gamma^{- 1} \left( \mathe^{n +
  \theta} \right) - \frac{1}{2} \right) \varphi_n = t_n \varphi_n
\end{equation}
where $t_n = t_n^{\ast}$ resp. $t_n^{\ast \ast}$. (in thats latter case with
$\theta$ lowered by 1/2)

Of course Eq.(\ref{GT}) for the operator $T$ is more appealing than (\ref{ts},
\ref{tss}) (where $n$ is replaced by $H$ and $t_n^{\ast}$ resp. $t_n^{\ast
\ast}$ are replaced by $T$) due to the combinatorial nature of the Gamma
function, but the eigenvalues of the operators are the same in the
``termodynamic limit'', $t \rightarrow \infty$.

\begin{figure}[h]
  \includegraphics[width=8cm]{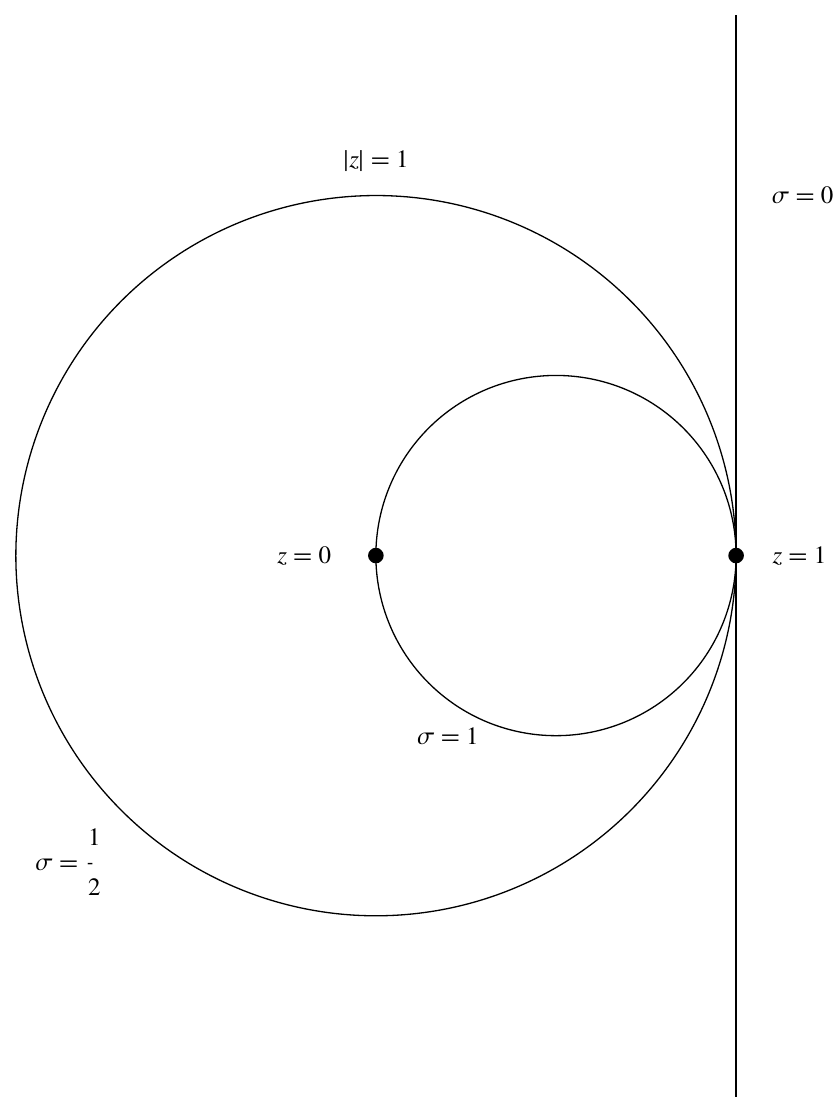}
  \caption{\label{accu}z-plane}
\end{figure}

{\tmstrong{Remark:}} If one consider the usual map $z \rightarrow 1 -
\frac{1}{s}$ then the critical line $s = \frac{1}{2} + i t$ ($t \in
\mathbbm{R}$) is mapped onto the unit circle $|z| = 1$; the two sets of
trivial zeros $\{t_n^{\ast} \}$ and $\{t_n^{\ast \ast} \} ^{}$ have as
accumulation point $z = 1$ (as $n \rightarrow \infty$), which is the same
accumulation point for the real zeros of the $\zeta$ function given by
$\bar{z}_n = 1 - \frac{1}{- 2 n} = 1 + \frac{1}{2 n}$, as $n \rightarrow
\infty$ (see Fig .\ref{accu}).

Neglecting the real zeros $\{z_n \}$, Fig. \ref{accu} illustrate by means of
two sets of trivial zeros $\{t_n^{\ast} \}$ and $\{t_n^{\ast \ast} \}$ the
Lee-Yang Theorem for the zeros of the partition function for some general spin
lattice system studied in statistical mechanics. If RH is true, then all non
trivial zeros of $\zeta (s)$ shall be located at the same circle \ $|z| = 1$,
with $z = 1$ as accumulation point.

\section{Conclusion}

In this note we have specified an operator equation for the operator $T$,
having as eigenvalues the trivial zeros of the $\zeta$ function on the
critical line (\ref{ts} and \ref{tss}); such zeros are the so called trivial
zeros (mean staircases), and are related in a strong non linear way to the
eigenvalues of a discrete harmonic oscillator described by $H$ in the large
$t$ limit. The introduction of $H$ may seem to be artificial: nevertheless the
two sets of trivial zeros ($\{t_n^{\ast} \}$ and $\{t_n^{\ast \ast} \}$) are
the eigenvalues of $T (H)$. In a subsequent note we will present a study of
another sequence of zeros possibly more connected with the true Riemann zeros.

\par
\vskip 0.3cm
\noindent
{\bfseries Update}
\par

Very recently G. Sierra and P.K. Townsend {\cite{sierra-2008}} introduced and
studied an interesting physical model (a charged particle in the plane in
presence of an electrical and a magnetic potential). In particular, the lowest
Landau level is connected with the smoothed counting function that gives the
average number of zeros, i.e. the staircase which here we have studied, by
means of a classical one-dimentional model of $N$ interacting charged
particles.

\par
\vskip 0.3cm
\noindent
{\bfseries Note added}
\par

Very recently Schumayer et al {\cite{schumayer-2008}} constructed (in
particular) the Quantum mechanical potential for $\xi (s)$ zeros, with the
first 200 energy eigenvalues (non trivial zeros). It is expected that the same
form of a quantum mechanical potential would appear using only the two sets of
zeroes we have discussed in this note. For the construction of an Hamiltonian
whose spectrum coincides with the primes, see also the recent work of S.
Sekatskii {\cite{sekatskii-2007}}.

\section{Appendix}

The hermitean matrix associated with the Mehta-Dyson model, $H$: discrete
annihilation and creation operators associated to $H$ whose spectrum is given
by the set of integers $(1, 2, \ldots, n)$, for any finite $n$.

In Ref{\cite{merlini-1999,bernasconi-2002}} it was studied the one dimentional
Mehta-Dyson model defined by the Hamiltonian $E = \overset{N}{\underset{i =
0}{\sum}} \left( \frac{1}{2} y_i^2 - \underset{i < j < N}{\sum} \log \left( |
y_i - y_j | \right) \right)$, where $y_i$ is the position of the i-th particle
on the line; then the fluctuation around the equilibrium positions (these are
given by the zeros of the Hermite polynomials of degree $N$, where $N$ is the
number of particles on the line, for every finite $N$), i.e. the harmonic
fluctuation spectrum is given by the eigenvalues of the hermitean $N \times N$
real matrix whose elements are given by:
\[  \left\{ \begin{array}{lll}
     H_{\tmop{ij}} = \frac{- 1}{| x_i - x_j |^2} &  & i \neq j\\
     &  & \\
     H_{\tmop{ij}} = 1 + \underset{k \neq i}{\sum} \frac{1}{| x_i - x_k |^2} &
     & i = j
   \end{array} \right. \]
$i, j = 1, \ldots, N$, where now the $x_i$ are the ``equilibrium positions''
i.e. the zeros of the Hermite polynomials of degree $N$.

The spectrum of $H$ is given exactly by the integers (1, 2, ... $N$) for every
finite $N$ and the eigenfunctions are given in terms of the Mehta-Dyson
polynomials of order 1 up to $N$. The Hamiltonian describing the harmonic
fluctuations takes then the form {\cite{merlini-1999}}:
\[ H = N \cdot I - \frac{1}{2} a a^{\ast} \]
where $I$ is the unit matrix of order $N$ and $a$, resp $a^{\ast}$, are the
discrete annihilation and creation operators (matrices of order $N \times N$)
which satisfy the commutator relation $[a, a^{\ast}] = - 2$.

Moreover $\left[ H, a^{\ast} \right] = a^{\ast}$ and $\left[ H, a \right] = -
a$.

If $X_k$ is the $k$-ten eigenvector of $H$ with eigenvalue the integer $k$,
one has:
\[ a^{\ast} X_{k + 1} = X_{k + 2} \]
and
\[ a X_{k + 1} = 2 \left( N - k \right) X_k \]
Explicitly, if $X_k = \left( \varphi_{1 k} \left( x_1 \right), \ldots,
\varphi_{\tmop{Nk}} \left( x_N \right) \right)$ is the $k$-ten eigenvector,
where $\varphi_k (x)$ is the k-th Mehta-Dyson polynomial of argument $x$, then
\[ a \varphi_{k + 1} \left( x_1 \right) = \frac{d}{d x_1}  \left( \varphi_{k +
   1} \left( x_1 \right) \right) = \underset{i \neq 1}{\overset{N}{\sum}}
   \frac{\varphi_{k + 1} \left( x_1 \right) - \varphi_{k + 1} \left( x_i
   \right)}{\left( x_1 - x_i \right)} \]
and
\[ a^{\ast} \varphi_{k + 1} \left( x_1 \right) = \left( 2 x_1 - \frac{d}{d
   x_1} \right) \varphi_{k + 1} \left( x_1 \right) \]
$a$ and $a^{\ast}$ are the two discrete annihilation and creation operators of
$H$.

$H$ as above with $N = n$ may be used to give the first $n$ trivial zeros of
the first set in (\ref{ts}) i.e. $t_1^{\ast} \ldots t_n^{\ast}$ while $H +
\frac{1}{2}$ may be used for obtaining the first $n$ trivial zeros of the
second set in (\ref{tss}) i.e. $t_1^{\ast \ast} \ldots t_n^{\ast \ast}$ in the
discussion on the mean staircases given in Section 4 above.

\end{document}